# Constrained Ordered Equilibrium Problems

Jinlu Li
Department of Mathematics
Shawnee State University
Portsmouth, Ohio 45662
USA

**Abstract**

In this paper, we consider some equilibrium problems (or saddle point problems), in which the domains of the considered mappings are limited at some regions. These restricted regions are defined by some mappings which are called the constrained mappings in the given problems. For this reason, we introduce the concepts of constrained ordered equilibrium problems, that are useful extensions of ordered equilibrium problems and ordinary equilibrium problems. Then, by using some fixed point theorems on posets, we prove several theorems for existence of solutions to some constrained ordered equilibrium problems. In particular, we investigate the solvability of some constrained ordered equilibrium problems in some partially ordered Banach spaces.

## 1. Introduction

Let $C$ and $D$ be nonempty sets. Let $T: C \times D \to R$ be a single-valued function. The ordinary equilibrium problem (or the saddle point problem) associated with the function $T$ and domains $C$, $D$, denoted by EP($T$, $C$, $D$), is to find $x^* \in C$, and $y^* \in D$ such that

$$T(x, y^*) \leq T(x^*, y^*) \leq T(x^*, y), \text{ for all } x \in C \text{ and } y \in D. \tag{1}$$

Where $(x^*, y^*) \in C \times D$ is called a solution to the problem EP($T$, $C$, $D$). The original results about the existence of solutions to EP($T$, $C$, $D$) was proved by Kakutani's fixed point theorem for set-valued mappings, in which $C$ and $D$ are supposed to be subsets of finite dimensional Euclidean spaces (see [12]).

In game theory, economic theory and some industry fields, the mapping $T$ in (1) is called a objective mapping. Sometimes the mapping $T$ may have incomplete utilities. That is, the values of the mapping $T$ may be in a nonlinear ordered poset. It leads authors to extend the ordinary equilibrium problem EP($T$, $C$, $D$) to ordered equilibrium problem as follows (see [5]). Let

$(U, \succcurlyeq^U)$ be a poset that is considered as the objective set of the considered problem. Let $T: C \times D \to R$ be a single-valued mapping. The ordered equilibrium problem (or the generalized saddle point problem) associated with the function $T$ and domains $C, D$, denoted by OEP($T, C, D$), is to find $x^* \in C$, and $y^* \in D$ such that

$$T(x, y^*) \not\succ^U T(x^*, y^*) \not\succ^U T(x^*, y), \text{ for all } x \in C \text{ and } y \in D. \qquad (2)$$

In particular, if $U$ is a topological vector space equipped with a partial order $\succcurlyeq^U$ that is induced by a closed and convex cone $K$ in $U$, then the ordered equilibrium problem OEP($T, C, D$) in (2) becomes a vector equilibrium problem VEP($T, C, D$) (see [4]): to find $x^* \in C$, and $y^* \in D$ such that

$$T(x, y^*) - T(x^*, y^*) \notin K/\{0\} \text{ and } T(x^*, y^*) - T(x^*, y) \notin K/\{0\}, \text{ for all } x \in C \text{ and } y \in D. \quad (3)$$

Similar to the ordinary equilibrium problems being an important branch in ordinary optimization theory (see [8], [12]), the vector equilibrium and ordered equilibrium problems have become a significant part in vector optimization theory and ordered equilibrium theory. These problems have received great attention from many authors (see [1-6], [10-11]).

Looking at the inequalities or order-inequalities in (1), (2), and (3), they are true for all $x \in C$ and $y \in D$, in which $C$ and $D$ are the underlying spaces in the given corresponding equilibrium problems. So, we can consider the problems (1), (2), and (3) as global equilibrium problems and the solutions $(x^*, y^*)$ are considered as global solutions. In this paper, we consider some local equilibrium problems, in which the corresponding inequalities are subjected to some given regions. These regions are related to the restricted mappings in the given problems. Such local equilibrium problems are called constrained ordered equilibrium problems (see section 3).

This paper is organized as follows: in section 2, we briefly recall some concepts of order-monotone mappings and a fixed point theorem on posets for easy reference; in section 3, we introduce the concepts of constrained ordered equilibrium problems and provide several existence theorems for their solutions; in section 4, we focus on constrained ordered equilibrium problems on partially ordered Banach spaces, which can be considered as special cases and more practicable in applications to game theory, economic theory, and optimization theory.

## 2. Preliminaries

Let $(Z, \succcurlyeq^Z)$, $(U, \succcurlyeq^U)$ be posets and $T: Z \to 2^U \setminus \{\emptyset\}$ a set-valued mapping. $T$ is said to be isotone, or to be $\succcurlyeq^Z$-increasing upward whenever, for any $x \preccurlyeq^Z y$ in $Z$ and $z \in Tx$, there is a $w \in Ty$ such that $z \preccurlyeq^U w$. $T$ is said to be $\succcurlyeq^Z$-increasing downward, if $x \preccurlyeq^Z y$ in $Z$, and $w \in Ty$, there is a $z \in Tx$ such that $z \preccurlyeq^U w$. If $T$ is both $\succcurlyeq^Z$-increasing upward and downward, then $T$ is said to be $\succcurlyeq^Z$-increasing.

$T$ is said to be $\succcurlyeq^Z$-decreasing upward whenever for any $x \preccurlyeq^Z y$ in $Z$ and $z \in Tx$, there is a $w \in Ty$ such that $z \succcurlyeq^U w$. $T$ is said to be $\succcurlyeq^Z$-decreasing downward, if $x \preccurlyeq^Z y$ in $Z$, and $w \in Ty$, then there is a $z \in Tx$ such that $z \succcurlyeq^U w$. If $T$ is both $\succcurlyeq^Z$-decreasing upward and downward, then $T$ is said to be $\succcurlyeq^Z$-decreasing.

Particularly, a single-valued mapping $F: Z \to U$ is said to be $\succcurlyeq^Z$-increasing whenever, for $x, y \in Z$, $x \preccurlyeq^Z y$ implies $F(x) \preccurlyeq^U F(y)$. An $\succcurlyeq^Z$-increasing mapping $F$ is said to be strictly $\succcurlyeq^Z$-increasing whenever $x \prec^Z y$ implies $F(x) \prec^U F(y)$. It can similarly be defined for $F$ to be $\succcurlyeq^Z$-decreasing and to be strictly $\succcurlyeq^Z$-decreasing.

In next section, by using Theorem 3.1 in [9], we prove an existence of solutions to some constrained ordered equilibrium problems and provide the inductive properties of the solution sets. We recall this theorem below for easy reference.

**Theorem 3.1** [9] *Let $(P, \succcurlyeq)$ be a chain-complete poset and let $F : P \to 2^P \setminus \{\emptyset\}$ be a set-valued mapping satisfying the following three conditions*:

A1. *F is $\succcurlyeq$-increasing upward.*
A2. *$(F(x), \succcurlyeq)$ is universally inductive, for every $x \in P$.*
A3. *There is an element $y'$ in $P$ and $v' \in F(y')$ with $y' \preccurlyeq v'$.*

*Then*

(i) *$(\mathcal{F}(F), \succcurlyeq)$ is a nonempty inductive poset;*

(ii) *$(\mathcal{F}(F) \cap [y'), \succcurlyeq)$ is a nonempty inductive poset; and $F$ has an $\succcurlyeq$-maximal fixed point $x'$ with $x' \succcurlyeq y'$.*

## 3. Constrained ordered equilibrium problems

### 3.1. Definition of constrained ordered equilibrium problems

Let $(X, \succcurlyeq^X)$, $(Y, \succcurlyeq^Y)$, and $(U, \succcurlyeq^U)$ be posets. Let $C$ and $D$ be nonempty subsets of $X$ and $Y$, respectively. Let $T: C \times D \to U$ be a single-valued mapping and $F: C \to 2^D \setminus \{\emptyset\}$, $G: D \to 2^C \setminus \{\emptyset\}$ be set-valued mappings. *The constrained ordered equilibrium problem* associated with domains $C$, $D$, utility set $U$, objective mapping $T$, and restriction mappings $F$, $G$, denoted by ROEP($T$, $F$, $G$), is to find $x^* \in C$, and $y^* \in D$ with

$$x^* \in G(y^*) \text{ and } y^* \in F(x^*) \text{ such that}$$

$$T(x, y^*) \not\succ^U T(x^*, y^*) \not\succ^U T(x^*, y), \text{ for all } x \in G(y^*) \text{ and } y \in F(x^*). \tag{4}$$

Where $(x^*, y^*) \in C \times D$ is called a solution to the problem ROEP($T$, $F$, $G$). Throughout this paper, unless otherwise is stated, we denote

$$S(T, F, G) = \text{the set of solutions to the problem ROEP}(T, F, G).$$

The following are special cases for problem ROEP($T$, $F$, $G$):

1. If both $F$ and $G$ are constant set-valued mappings with values $D$ and $C$, respectively, then the constrained ordered equilibrium problem (4) becomes an ordered equilibrium problem OEP($T$, $C$, $D$): to find $x^* \in C$, and $y^* \in D$, such that

$$T(x, y^*) \not\succ^U T(x^*, y^*) \not\succ^U T(x^*, y), \text{ for all } x \in C \text{ and } y \in D.$$

2. If $F$ is a constant set-valued mapping with value $D$ and $G: D \to 2^C \setminus \{\emptyset\}$, then the constrained ordered equilibrium problem (4) becomes a one side constrained ordered equilibrium problem ROEP($T$, $G$, $C$): to find $x^* \in C$, and $y^* \in D$ with

$x^* \in G(y^*)$ and $y^* \in D$ such that

$$T(x, y^*) \not\succ^U T(x^*, y^*) \not\succ^U T(x^*, y), \text{ for all } x \in G(y^*) \text{ and } y \in D.$$

3. If $G$ is a constant set-valued mapping with value $C$ and $F: C \to 2^D\setminus\{\varnothing\}$, then the constrained ordered equilibrium problem (4) becomes a one side constrained ordered equilibrium problem ROEP($T, F, D$): to find $x^* \in C$, and $y^* \in D$ with

$$x^* \in C \text{ and } y^* \in F(x^*) \text{ such that}$$

$$T(x, y^*) \not\succ^U T(x^*, y^*) \not\succ^U T(x^*, y), \text{ for all } x \in C \text{ and } y \in F(x^*).$$

4. If $(U, \succ^U) = (R, \geq)$, then the constrained ordered equilibrium problem (4) becomes a constrained (normal) equilibrium problem REP($T, F, G$): to find $x^* \in C$, and $y^* \in D$, with

$$x^* \in G(y^*) \text{ and } y^* \in F(x^*) \text{ such that}$$

$$T(x, y^*) \leq T(x^*, y^*) \leq T(x^*, y), \text{ for all } x \in G(y^*) \text{ and } y \in F(x^*).$$

It is equivalent to

$$\sup_{x \in G(y^*)} T(x, y^*) = T(x^*, y^*) = \inf_{y \in F(x^*)} T(x^*, y).$$

As in the above case 2 and case 3, we can consider one side constrained equilibrium problems. In particular, we have

5. If both $F$ and $G$ are constant set-valued mappings with values $D$ and $C$, respectively, and $(U, \succ^U) = (R, \geq)$, then the constrained ordered equilibrium problem (1) becomes an (normal) equilibrium problem EP($T, C, D$): to find $x^* \in C$, and $y^* \in D$ such that

$$T(x, y^*) \leq T(x^*, y^*) \leq T(x^*, y), \text{ for all } x \in C \text{ and } y \in D.$$

It is equivalent to

$$\sup_{x \in C} T(x, y^*) = T(x^*, y^*) = \inf_{y \in D} T(x^*, y).$$

For a given objective mapping $T: C \times D \to U$, and constrained mappings $F: C \to 2^D\setminus\{\varnothing\}$, and $G: D \to 2^C\setminus\{\varnothing\}$, we define their order-optimization mappings $\varphi: C \to 2^D\setminus\{\varnothing\}$ and $\psi: D \to 2^C\setminus\{\varnothing\}$ as follows:

$$\varphi(x) = \{yz \in F(x): T(x, yz) \text{ is an } \succ^U\text{-minimal point in } T(x, F(x))\}, \text{ for every } x \in C;$$
$$\psi(y) = \{xu \in G(y): T(xu, y) \text{ is an } \succ^U\text{-maximal point in } T(G(y), y)\}, \text{ for every } y \in D.$$

In particular, if both $F$ and $G$ are constant set-valued mappings with values $D$ and $C$, respectively, then the order-optimization mappings for an objective mapping $T$ is defined as follows:

$$\Phi(x) = \{z \in D: T(x, z) \text{ is an } \succ^U\text{-minimal point in } T(x, D)\}, \text{ for every } x \in C;$$
$$\Psi(y) = \{u \in C: T(u, y) \text{ is an } \succ^U\text{-maximal point in } T(C, y)\}, \text{ for every } y \in D.$$

### 3.2. Solvability of some constrained ordered equilibrium problems

In this subsection, unless otherwise is stated, we let $(X, \succcurlyeq^X)$ and $(Y, \succcurlyeq^Y)$ be posets and let $C$ and $D$ be chain-complete subsets in $(X, \succcurlyeq^X)$ and $(Y, \succcurlyeq^Y)$, respectively. Let $T: C \times D \to U$ be an objective mapping and let $F: C \to 2^D \setminus \{\emptyset\}$ and $G: D \to 2^C \setminus \{\emptyset\}$ be constrained mappings. Let $\varphi: C \to 2^D \setminus \{\emptyset\}$ and $\psi: D \to 2^C \setminus \{\emptyset\}$ be their order-optimization mappings.

**Theorem 3.1.** *For the mappings $T$, $F$, and $G$, suppose that the order-optimization mappings $\varphi$ and $\psi$ both are order-increasing upward with universally inductive values in $(D, \succcurlyeq^Y)$ and $(C, \succcurlyeq^X)$, respectively. If there are $(x', y') \in C \times D$, values $u' \in \varphi(x')$ and $z' \in \psi(y')$ such that*

$$x' \preccurlyeq^X z' \quad \text{and} \quad y' \preccurlyeq^Y u', \tag{5}$$

*then, $S(T, F, G)$ is a nonempty and inductive subset of $(C \times D, \succcurlyeq^{X \times Y})$, and $\mathrm{ROEP}(T, F, G)$ has an $\succcurlyeq^{X \times Y}$-maximal solution.*

*Proof.* Let $\succcurlyeq^{X \times Y} = \succcurlyeq^X \times \succcurlyeq^Y$ be the component-wise ordering on $C \times D$ satisfying that, for $(x_1, y_1)$, $(x_2, y_2) \in C \times D$,

$$(x_1, y_1) \preccurlyeq^{X \times Y} (x_2, y_2) \text{ if and only if, } x_1 \preccurlyeq^X x_2 \text{ and } y_1 \preccurlyeq^Y y_2.$$

Since $(C, \succcurlyeq^X)$ and $(D, \succcurlyeq^Y)$ both are chain-complete subsets in $(X, \succcurlyeq^X)$ and $(Y, \succcurlyeq^Y)$, respectively, then $(C \times D, \succcurlyeq^{X \times Y})$ is an $\succcurlyeq^{X \times Y}$-chain-complete poset. By using the order-optimization mappings $\varphi$ and $\psi$, we define $\Gamma: C \times D \to 2^{C \times D}$ by

$$\Gamma(x, y) = \psi(y) \times \varphi(x), \text{ for } (x, y) \in C \times D.$$

From the conditions about the order-optimization mappings in this theorem, it implies that $\Gamma: C \times D \to 2^{C \times D} \setminus \{\emptyset\}$ is well defined.

We first show that $\Gamma: C \times D \to 2^{C \times D} \setminus \{\emptyset\}$ is $\succcurlyeq^{X \times Y}$-increasing upward. For any $(x_1, y_1)$, $(x_2, y_2) \in C \times D$ with $(x_1, y_1) \preccurlyeq^{X \times Y} (x_2, y_2)$, which is equivalent to both $x_1 \preccurlyeq^X x_2$ and $y_1 \preccurlyeq^Y y_2$. For an arbitrary $(z_1, u_1) \in \Gamma(x_1, y_1)$, it is equivalent to both:

$$z_1 \in \psi(y_1) \subseteq C \text{ and } u_1 \in \varphi(x_1) \subseteq F(x_1) \subseteq D.$$

Since $\psi: D \to 2^C \setminus \{\emptyset\}$ is $\succcurlyeq^Y$-increasing upward, from $y_1 \preccurlyeq^Y y_2$ and $z_1 \in \psi(y_1)$, there is $z_2 \in \psi(y_2)$ such that

$$z_1 \preccurlyeq^X z_2. \tag{6}$$

Since $\varphi: C \to 2^D \setminus \{\emptyset\}$ is $\succcurlyeq^X$-increasing upward, from $x_1 \preccurlyeq^X x_2$ and $u_1 \in \varphi(x_1)$, there is $u_2 \in \varphi(x_2)$ such that

$$u_1 \preccurlyeq^Y u_2. \tag{7}$$

From (7) and (6), we get

$$(z_2, u_2) \in \psi(y_2) \times \varphi(x_2) = \Gamma(x_2, y_2) \text{ such that } (z_1, u_1) \preccurlyeq^{X \times Y} (z_2, u_2).$$

Hence $\Gamma: C \times D \to 2^{C \times D} \setminus \{\emptyset\}$ is $\succcurlyeq^{X \times Y}$-increasing upward, that satisfies condition (A1) in Theorem 3.1 in [15].

Then, we prove that, for every $(x, y) \in C \times D$, $(\Gamma(x, y), \succcurlyeq^{X \times Y})$, there is an universally inductive sub-poset of $(C \times D, \succcurlyeq^{X \times Y})$. For an arbitrary $\succcurlyeq^{X \times Y}$-chain $\{(x_\alpha, y_\alpha)\} \subseteq (C \times D, \succcurlyeq^{X \times Y})$ satisfying that, for every point $(x_\alpha, y_\alpha)$ in this chain, there is $(z_\alpha, u_\alpha) \in \Gamma(x, y) = \psi(y) \times \varphi(x)$ such that

$$(x_\alpha, y_\alpha) \preccurlyeq^{X \times Y} (z_\alpha, u_\alpha).$$

It is equivalent to both of the following

$$x_\alpha \preccurlyeq^X z_\alpha \in \psi(y) \text{ and } y_\alpha \preccurlyeq^Y u_\alpha \in \varphi(x). \tag{8}$$

Since $\{x_\alpha\}$ is an $\succcurlyeq^X$-chain in $(C, \succcurlyeq^X)$, from (8) and conditions about the order-optimization mapping $\psi$ in this theorem, $(\psi(y), \succcurlyeq^X)$ is universally inductive in $(C, \succcurlyeq^X)$. There is $z \in \psi(y)$ which is an $\succcurlyeq^X$-upper bound of the chain $\{y_\alpha\}$. That is,

$$x_\alpha \preccurlyeq^X z \in \psi(y), \text{ for all } \alpha. \tag{9}$$

Since $\{y_\alpha\}$ is an $\succcurlyeq^Y$-chain in $(D, \succcurlyeq^Y)$, from (8) and condition about the order-optimization mapping $\varphi$ in this theorem, $(\varphi(x), \succcurlyeq^Y)$ is universally inductive in $(D, \succcurlyeq^Y)$, there is $u \in \varphi(x)$ which is an $\succcurlyeq^Y$-upper bound of the chain $\{y_\alpha\}$. That is,

$$y_\alpha \preccurlyeq^X u \in \varphi(x), \text{ for all } \alpha. \tag{10}$$

From (10) and (9), the element $(z, u) \in \Gamma(x, y)$ satisfies

$$(x_\alpha, y_\alpha) \preccurlyeq^{X \times Y} (z, u), \text{ for all } \alpha.$$

It implies that $(\Gamma(x, y), \succcurlyeq^{X \times Y})$ is universally inductive in $(C \times D, \succcurlyeq^{X \times Y})$. Hence, $\Gamma$ satisfies condition (A2) in Theorem 3.1 in [15].

For the given point $(x', y')$ in (5) in this theorem, it implies

$$(z', u') \in \psi(y') \times \varphi(x') = \Gamma(x', y') \text{ with } (x', y') \preccurlyeq^{X \times Y} (z', u').$$

Hence $\Gamma$ satisfies condition (A3) in Theorem 3.1 in [15], from which it implies that $(\mathcal{F}(\Gamma), \succcurlyeq^{X \times Y})$ is a nonempty inductive poset contained in $(C \times D, \succcurlyeq^{X \times Y})$.

For any fixed point $(x^*, y^*) \in \mathcal{F}(\Gamma)$, we have

$$(x^*, y^*) \in \Gamma(x^*, y^*) = \psi(y^*) \times \varphi(x^*).$$

That is,

$$x^* \in \psi(y^*) = \{xz \in G(y^*): T(xz, y^*) \text{ is an } \succcurlyeq^U\text{-maximal point in } T(G(y^*), y^*)\};$$

and $\quad y^* \in \varphi(x^*) = \{yu \in F(x^*): T(x^*, yu) \text{ is an } \succcurlyeq^U\text{-minimal point in } T(x^*, F(x^*))\}.$

It follows that $x^* \in G(y^*)$ satisfying

$$T(x^*, y^*) \text{ is an } \succcurlyeq^U\text{-maximal point in } T(G(y^*), y^*) = \{T(x, y^*): x \in G(y^*)\};$$

and $y^* \in F(x^*)$ satisfying

$$T(x^*, y^*) \text{ is an } \geqslant^U\text{-minimal point in } T(x^*, F(x^*)) = \{T(x^*, y): y \in F(x^*)\}.$$

It implies that the point $(x^*, y^*)$ satisfies (1) and it is a solution to ROEP($T, F, G$). Hence, the solution set of ROEP($T, F, G$) is a nonempty inductive sub-poset of $(C \times D, \geqslant^{X \times Y})$. □

As a consequence of Theorem 3.1, we have the following results.

**Corollary 3.2**. *For the mappings $T, F$ and $G$, suppose that the order-optimization mappings $\varphi$ and $\psi$ both are order-increasing with universally bi-inductive values in $(D, \geqslant^Y), (C, \geqslant^X)$, respectively. If there are points $(x', y'), (x'', y'') \in C \times D$, with $(x', y') \leqslant^{X \times Y} (x'', y'')$, and values $u' \in \varphi(x'), u'' \in \varphi(x''), z' \in \psi(y')$, and $z'' \in \psi(y'')$ satisfying (5) and*

$$x'' \geqslant^X z'' \quad \text{and} \quad y'' \geqslant^Y u'',$$

*then, $\mathcal{S}(T, F, G)$ is nonempty and bi-inductive in $(C \times D, \geqslant^{X \times Y})$, and ROEP($T, F, G$) has both $\geqslant^{X \times Y}$-maximal and $\geqslant^{X \times Y}$-minimal solutions.*

Conditions in Corollary 3.2 can be substituted by some conditions on $C$ and $D$.

**Corollary 3.3**. *For the mappings $T, F$ and $G$, suppose that the order-optimization mappings $\varphi$ and $\psi$ both are order-increasing with universally bi-inductive values in $(D, \geqslant^Y)$ and $(C, \geqslant^X)$, respectively. Suppose that both $C$ and $D$ possess order-maximum and order-minimum points. Then, $\mathcal{S}(T, F, G)$ is nonempty and bi-inductive in $(C \times D, \geqslant^{X \times Y})$, and ROEP($T, F, G$) has both $\geqslant^{X \times Y}$-maximal and $\geqslant^{X \times Y}$-minimal solutions.*

In particular, in the order-optimization mappings $\varphi$ and $\psi$ both are single valued, since singletons are automatically considered as *universally inductive*, we obtain the following result as an immediate consequence of Theorem 3.1.

**Corollary 3.4**. *For the mappings $T, F$ and $G$, suppose that the order-optimization mappings $\varphi$ and $\psi$ both are order-increasing single-valued mappings. If there is $(x', y') \in C \times D$ such that*

$$x' \leqslant^X \psi(y') \quad \text{and} \quad y' \leqslant^Y \varphi(x'),$$

*then, $\mathcal{S}(T, F, G)$ is nonempty and inductive in $(C \times D, \geqslant^{X \times Y})$, and ROEP($T, F, G$) has an $\geqslant^{X \times Y}$-maximal solution.*

**Corollary 3.5**. *Suppose that the order-optimization mappings of an objective mapping $T$, $\Phi$ and $\Psi$ both are order-increasing upward with universally inductive values in $(D, \geqslant^Y)$ and $(C, \geqslant^X)$, respectively. If there are $(x', y') \in C \times D$, values $u' \in \Phi(x')$ and $z' \in \Psi(y')$ satisfying (5), then, the problem OEP($T, C, D$) has an $\geqslant^{X \times Y}$-maximal solution. Moreover, the solution set is inductive in $(C \times D, \geqslant^{X \times Y})$.*

## 4. Solvability of constrained ordered equilibrium problems on partially ordered Banach spaces

In this section, we consider the solvability of constrained ordered equilibrium problems on partially ordered Banach spaces, which can be considered as special cases of the problems studied in the previous sub-section.

Let $(U, \succcurlyeq^U)$ be a poset. In this subsection, unless otherwise is stated, we let $(X, \succcurlyeq^X)$ and $(Y, \succcurlyeq^Y)$ be partially ordered Banach spaces, in which the partial orders $\succcurlyeq^X$ and $\succcurlyeq^Y$ are induced by closed and convex cones $K_X$ and $K_Y$ in $X$ and $Y$, respectively and let $C$ and $D$ be nonempty subsets of $X$ and $Y$, respectively. Then, with an objective mapping $T: C \times D \to U$ and constrained mappings $F: C \to 2^D \setminus \{\emptyset\}$ and $G: D \to 2^C \setminus \{\emptyset\}$, we consider the constrained ordered equilibrium problem ROEP($T, F, G$).

In Theorem 3.1 in [9], the underlying poset $(P, \succcurlyeq)$ is required to be $\succcurlyeq$-chain-complete. As we consider partially ordered Banach spaces as underlying spaces, some results regarding to chain-complete and universally inductive from [6] and [9] are useful, which are summarized as the following lemma.

**Lemma 4.1.** *Every nonempty weakly compact subset of a partially ordered Banach space is chain-complete and universally inductive.*

**Theorem 4.2**. *Let $C$ and $D$ be weakly compact subsets in partially ordered Banach spaces $(X, \succcurlyeq^X)$ and $(Y, \succcurlyeq^Y)$, respectively. For the mappings $T$, $F$ and $G$, suppose that the order-optimization mappings $\varphi$ and $\psi$ both are order-increasing upward with closed values in $D$ and $C$, respectively. If there are $(x', y') \in C \times D$, values $u' \in \varphi(x')$ and $z' \in \psi(y')$ satisfying (5), then, $S(T, F, G)$ is nonempty and inductive in $(C \times D, \succcurlyeq^{X \times Y})$, and ROEP($T, F, G$) has an $\succcurlyeq^{X \times Y}$-maximal solution.*

*Proof.* Notice that a subset in a Banach space is closed with respect to the norm topology, if and only if, it is closed with respect to the weak topology. It follows that the norm topology is natural with respect to a given partial order on this Banach space if and only if the weak topology is nature with respect to the same partial order. Then, from Lemma 4.1, $(C, \succcurlyeq^X)$ and $(D, \succcurlyeq^Y)$ are chain-complete, so is $(C \times D, \succcurlyeq^{X \times Y})$. Since $\varphi$ and $\psi$ both have closed values in $D$ and $C$ which are weakly compact subsets in $X$ and $Y$, respectively, it implies that the values of $\varphi$ and $\psi$ are also weakly compact. By Lemma 4.1, they are universally inductive. Then this theorem follows from Theorem 3.1 immediately. □

Since every closed, convex, and (norm) bounded subset in a reflexive Banach space is weakly compact, then the following corollary immediately follows from Theorem 4.2.

**Corollary 4.3**. *Let $C$ and $D$ be nonempty closed convex and (norm) bounded subsets in partially ordered reflexive Banach spaces $(X, \succcurlyeq^X)$ and $(Y, \succcurlyeq^Y)$, respectively. Suppose that the order-optimization mappings $\varphi$ and $\psi$ of mappings $T$, $F$ and $G$ both are order-increasing upward with closed values in $D$ and $C$, respectively. If there are $(x', y') \in C \times D$, values $u' \in \varphi(x')$ and $z' \in \psi(y')$ satisfying (5), then, $S(T, F, G)$ is nonempty and inductive in $(C \times D, \succcurlyeq^{X \times Y})$; ROEP($T, F, G$) has an $\succcurlyeq^{X \times Y}$-maximal solution.*

To conclude this section, we consider constrained equilibrium problems on partially ordered Banach spaces, in which $T: C \times D \to R$ is a real valued function. We have the following corollaries of Theorem 4.2.

**Corollary 4.4**. *Let C and D be weakly compact subsets in partially ordered Banach spaces $(X, \succcurlyeq^X)$ and $(Y, \succcurlyeq^Y)$, respectively. Let $T: C \times D \to R$ be a real valued function. For the mappings T, F and G, suppose that the order-optimization mappings $\varphi$ and $\psi$ both are order-increasing upward with closed values in D and C, respectively. If there are $(x', y') \in C \times D$, values $u' \in \varphi(x')$ and $z' \in \psi(y')$ such that (11.1) holds, then, the problem REP(T, F, G) has a solution. That is, there are $x^* \in C$, and $y^* \in D$, with $x^* \in G(y^*)$ and $y^* \in F(x^*)$ such that*

$$\sup_{x \in G(y^*)} T(x, y^*) = T(x^*, y^*) = \inf_{y \in F(x^*)} T(x^*, y).$$

*Moreover, the solution set of REP(T, F, G) is inductive in $(C \times D, \succcurlyeq^{X \times Y})$.*

**Corollary 4.5**. *Let C and D be weakly compact subsets in partially ordered Banach spaces $(X, \succcurlyeq^X)$ and $(Y, \succcurlyeq^Y)$, respectively. Let $T: C \times D \to R$ be a real valued objective function. Suppose that its order-optimization mappings $\Phi$ and $\Psi$ both are order-increasing with closed values in D and C, respectively. If there are $(x', y') \in C \times D$, values $u' \in \Phi(x')$ and $z' \in \Psi(y')$ such that $x' \preccurlyeq^X z'$ and $y' \preccurlyeq^Y u'$, then, , there are $x^* \in C$, $y^* \in D$, such that*

$$\sup_{x \in C} T(x, y^*) = T(x^*, y^*) = \inf_{y \in D} T(x^*, y).$$

*Moreover, the solution set of EP(T, C, D) is inductive in $(C \times D, \succcurlyeq^{X \times Y})$.*

## 5. Applications to constrained two-person zero sums strategic games

In this section, for an arbitrary positive integer $k$, let $(R^k, \succcurlyeq^k)$ be the $k$-dimensional Euclidean space, where $\succcurlyeq^k$ is the component-wise ordering on $R^k$. We consider a two-person strategic game with complete preferences. Suppose that players 1 and 2 have strategic sets (profile sets) C and D, which are nonempty closed and convex subsets of $R^m$ and $R^n$, respectively. Let $T: C \times D \to R$ be a real valued function that defines the utility of player 1. Then player 2 has utility function $-T$. It implies that this game is a zero sum game. Suppose that this game is played under the following constrained conditions which are defined by mappings $F: C \to 2^D \setminus \{\varnothing\}$ and $G: D \to 2^C \setminus \{\varnothing\}$ both having closed values. It means that, for every strategy $x \in C$ for player 1 to play, player 2 must be restricted to choose his strategy from $F(x) \subseteq D$ and, for every strategy $y \in D$ for player 2 to play, player 1 must be restricted to choose his strategy from $G(y) \subseteq C$. This game is called a constrained two-person zero sums strategic game and denoted by G(T, C, D, F, G).

As a consequence of Corollary 4.5, we have an existence of constrained equilibrium for some constrained games.

**Theorem 5.1**. *In a game G(T, C, D, F, G), suppose that C and D are bounded and the order-optimization mappings $\varphi$ and $\psi$ both are order-increasing upward with closed values in D and C, respectively. Suppose that there are $(x', y') \in C \times D$, values $u' \in \varphi(x')$ and $z' \in \psi(y')$ satisfying (5), then, this game has a constrained equilibrium $(x^*, y^*) \in C \times D$, with $x^* \in G(y^*)$ and $y^* \in F(x^*)$ such that*

$$T(x, y^*) \leq T(x^*, y^*) \leq T(x^*, y), \text{ for all } x \in G(y^*) \text{ and } y \in F(x^*).$$

*Moreover, the collection of constrained equilibria of this game is inductive in* $(C \times D, \succcurlyeq^{m+n})$.

**Observation 5.2**. It is clear that the conditions in Theorem 5.1 that the order-optimization mappings $\varphi$ and $\psi$ both have closed values in $D$ and $C$, respectively, can be obtained by the continuity of the utility function $T: C \times D \to R$.